\documentclass[12pt,letterpaper]{amsart}
\usepackage{amsmath,txfonts}
\usepackage{amssymb}
\usepackage{amsxtra}
\usepackage{amsthm, color}
\usepackage{txfonts}
\usepackage{graphicx}
\usepackage{times}
\usepackage{citeref}
\usepackage{tikz}
\usepackage{hyperref}
\usepackage{stmaryrd}
\usepackage[T3,T1]{fontenc}
\usepackage{pgfplots}

\usetikzlibrary{calc}

\usepackage{mathrsfs}
\usepackage{amsfonts}
\usepackage{amssymb}
\usepackage{ifthen}
\usepackage{graphicx}
\nonstopmode \numberwithin{equation}{section}

\newtheorem{thm}{Theorem}[section]
\newtheorem{lem}{Lemma}[section]%[thm]
\newtheorem{cor}[thm]{Corollary}
\newtheorem{prop}[thm]{Proposition}

\newtheorem{step}{Step}[section]%[thm]

\theoremstyle{definition}
\newtheorem{mlem}{Main lemma}[section]
\newtheorem{assertion}{Assertion}[section]
\newtheorem{cl}{Claim}[section]
\newtheorem{ca}{Case}[section]
\newtheorem{sca}{Subcase}[section]
\newtheorem{scl}{Subclaim}[section]
\newtheorem{conj}[thm]{Conjecture}
\newtheorem{fact}{Fact}[section]
\newtheorem{defn}[thm]{Definition}
\newtheorem{op}[thm]{Open Problem}
\newtheorem{prob}{Problem}[section]
\newtheorem{ques}{Question}[section]
\newtheorem{rem}[thm]{Remark}
\newtheorem{exam}[thm]{Example}

\numberwithin{equation}{section}

\newcounter {own}
\def\theown {\thesection       .\arabic{own}}

\newenvironment{pf}[1][]{%
 \vskip 3mm
 \noindent
 \ifthenelse{\equal{#1}{}}%
  {{\slshape Proof. }}%
  {{\slshape #1.} }%
 }%
{\qed\bigskip}

\newcounter{alphabet}
\renewcommand{\thealphabet}{\Alph{alphabet}}

\newenvironment{Thm}[1][]{\refstepcounter{alphabet}%
	\bigskip
	\noindent
	{\bf Theorem \thealphabet}%
	\ifthenelse{\equal{#1}{}}{}{ (#1)}%
	{\bf .} \itshape
}{\vskip 8pt}

\newenvironment{Lem}[1][]{\refstepcounter{alphabet}%
	\bigskip
	\noindent
	{\bf Lemma \thealphabet}%
	{\bf .} \itshape
}{\vskip 8pt}

\def\be{\begin{equation}}
\def\ee{\end{equation}}

\newcommand{\ben}{\begin{enumerate}}
\newcommand{\een}{\end{enumerate}}

\newcommand{\blem}{\begin{lem}}
\newcommand{\elem}{\end{lem}}
\newcommand{\bthm}{\begin{thm}}
\newcommand{\ethm}{\end{thm}}
\newcommand{\bcor}{\begin{cor}}
\newcommand{\ecor}{\end{cor}}
\newcommand{\beg}{\begin{exam}}
\newcommand{\eeg}{\end{exam}}
\newcommand{\begs}{\begin{examples}}
\newcommand{\eegs}{\end{examples}}
\newcommand{\bdefe}{\begin{defn}}
\newcommand{\edefe}{\end{defn}}
\newcommand{\bprob}{\begin{prob}}
\newcommand{\eprob}{\end{prob}}
\newcommand{\bques}{\begin{ques}}
\newcommand{\eques}{\end{ques}}
\newcommand{\bei}{\begin{itemize}}
\newcommand{\eei}{\end{itemize}}
\newcommand{\bcon}{\begin{conj}}
\newcommand{\econ}{\end{conj}}
\newcommand{\bop}{\begin{op}}
\newcommand{\eop}{\end{op}}

\newcommand{\bas}{\begin{assertion}}
\newcommand{\eas}{\end{assertion}}

\newcommand{\bfa}{\begin{fact}}
\newcommand{\efa}{\end{fact}}

\newcommand{\bca}{\begin{ca}}
\newcommand{\eca}{\end{ca}}

\newcommand{\bst}{\begin{step}}
\newcommand{\est}{\end{step}}

\newcommand{\bsca}{\begin{sca}}
\newcommand{\esca}{\end{sca}}

\newcommand{\bcl}{\begin{cl}}
\newcommand{\ecl}{\end{cl}}

\newcommand{\bmlem}{\begin{mlem}}
\newcommand{\emlem}{\end{mlem}}

\newcommand{\bscl}{\begin{scl}}
\newcommand{\escl}{\end{scl}}

\newcommand{\bcons}{\begin{conjs}}
\newcommand{\econs}{\end{conjs}}

\newcommand{\bprop}{\begin{prop}}
\newcommand{\eprop}{\end{prop}}

\newcommand{\br}{\begin{rem}}
\newcommand{\er}{\end{rem}}
\newcommand{\brs}{\begin{rems}}
\newcommand{\ers}{\end{rems}}
\newcommand{\bo}{\begin{obser}}
\newcommand{\eo}{\end{obser}}
\newcommand{\bos}{\begin{obsers}}
\newcommand{\eos}{\end{obsers}}
\newcommand{\bpf}{\begin{pf}}
\newcommand{\epf}{\end{pf}}
\newcommand{\ba}{\begin{array}}
\newcommand{\ea}{\end{array}}
\newcommand{\beq}{\begin{eqnarray}}
\newcommand{\beqq}{\begin{eqnarray*}}
\newcommand{\eeq}{\end{eqnarray}}
\newcommand{\eeqq}{\end{eqnarray*}}

\newcounter{minutes}
\divide\time by 60
\newcounter{hours}
\multiply\time by 60 \addtocounter{minutes}{-\time}
\begin{document}

%%%%%%%%%%%%%%%%%%%%%%%%%%
%%%%%%%%%%%%%%%%%%%%%%%%%%
%%% Added by M Vuorinen:
%\def\thefootnote{}
%\footnotetext{ \texttt{\tiny File:~\jobname .tex,
 %          printed: \number\year-\number\month-\number\day,
 %          \thehours.\ifnum\theminutes<10{0}\fi\theminutes}
%} \makeatletter\def\thefootnote{\@arabic\c@footnote}\makeatother
%%%%%%%%%%%%%%%%%%%%%%%%%%
%%%%%%%%%%%%%%%%%%%%%%%%%%

\bibliographystyle{amsplain}
\title []
{Sharp constant problems for harmonic mappings}

\author{Shaolin Chen$^{\ast}$}
\address{S. L. Chen,    Center for Applied Mathematics of Guangxi, College of Mathematics and Statistics, Guangxi Normal University,
Guilin, Guangxi 541004, People's Republic of China} \email{mathechen@126.com}

%\author{Manzi Huang}
%\address{M. Z.  Huang, Key Laboratory of High Performance Computing and Stochastic Information Processing,
%College of Mathematics and Statistics, Hunan Normal University, Changsha, Hunan 410081, People's Republic of China}
%\email{mzhuang@hunnu.edu.cn}

\author{Caicai Hu}
\address{C. C. Hu, School of Computational Science and Electronics, Hunan Institute of Engineering, Xiangtan, Hunan 411104, P. R. China}
\email{caicaihu33@163.com}

\author{Lei Jin}
\address{L. Jin, Key Laboratory of Computing and Stochastic Mathematics (Ministry of Education), 
College of Mathematics and Statistics, Hunan Normal University, Changsha, Hunan 410081, P. R. China}
\email{leijin0825@163.com}

\author{Qianyun Li}
	\address{Q. Y. Li, Key Laboratory of Computing and Stochastic Mathematics (Ministry of Education), 
College of Mathematics and Statistics, Hunan Normal University, Changsha, Hunan 410081, P. R. China}
	\email{liqianyun@hunnu.edu.cn}

\thanks{\\
${}^{\mathbf{*}}$ Corresponding author}

\begin{abstract}
 The primary objective of this paper is to investigate sharp constant problems in harmonic Hardy spaces. We first establish a sharp asymptotic constant for a Littlewood-Paley type inequality applicable to weak harmonic $(K, K_0)$-quasiregular mappings, where $K \geq 1$ and $K_0 \geq 0$. Subsequently, we derive sharp Hardy-Littlewood type constants for harmonic mappings.
\end{abstract}

\maketitle \pagestyle{myheadings} \markboth{S. L. Chen, C. C. Hu, L. Jin and Q. Y. Li}{Sharp constant problems for harmonic mappings}

\dedicatory{}

%%%%%%%%%%%%%%%%%%%%%%%%%%%%%%%%%%%%%%%%%%%%%%%%%%%%

\section{Introduction and main results}\label{sec-1}
Let $\mathbb{D} = \{ z : |z| < 1 \}$ denote the open unit disk in the complex plane $\mathbb{C}$, and let $\mathbb{T} := \partial \mathbb{D}$ be its boundary, the unit circle.  
For $z=x+iy\in\mathbb{C}$, the complex formal derivatives operators are defined by
$$\frac{\partial}{\partial z}=\frac{1}{2}\left(\frac{\partial }{\partial x}-i\frac{\partial }{\partial y}\right)
~~~\mbox{and}~~~\frac{\partial}{\partial \overline{z}}=\frac{1}{2}\left(\frac{\partial f}{\partial x}+i\frac{\partial f}{\partial y}\right).$$

For $z = r e^{i\vartheta} \in \mathbb{C}$ and $\vartheta \in [0, 2\pi]$, the directional derivative of $f$ at $z$ in the direction $\vartheta$ is defined as
$$
\partial_\vartheta f(z)
=
\lim_{r \to 0^+} \frac{f(z + r e^{i\vartheta}) - f(z)}{r}
=
e^{i\vartheta} f_z(z) + e^{-i\vartheta} f_{\bar z}(z),
$$ where $f_z=\partial f/\partial z$ and $f_{\overline{z}}=\partial f/\partial \overline{z}$.
Consequently, the maximum and minimum moduli of the directional derivative are given respectively by
$$
\Lambda_f(z) := \max_{0 \leq \vartheta \leq 2\pi} |\partial_\vartheta f(z)|
= |f_z(z)| + |f_{\bar z}(z)|,
$$
and
$$
\lambda_f(z) := \min_{0 \leq \vartheta\leq 2\pi} |\partial_\vartheta f(z)|
= \bigl| |f_z(z)| - |f_{\bar z}(z)| \bigr|.
$$

%\subsection*{Harmonic functions}

A twice continuously differentiable complex-valued function $f = u + iv$ is called a \emph{harmonic mapping} in the domain $ \Omega \subset \mathbb{C} $ if both its real and imaginary parts are harmonic, that is,
$$
\Delta u = \Delta v = 0 \quad \text{in } \Omega,
$$
where \(\Delta\) denotes the Laplacian operator
$$
\Delta := 4\frac{\partial^{2}}{\partial z\,\partial \overline{z}}
      = \frac{\partial^{2}}{\partial x^{2}} + \frac{\partial^{2}}{\partial y^{2}},
$$
and \(z = x + iy \in \Omega\). It is well-known that every harmonic
mapping $f$ defined in a simply connected domain $\Omega$ admits a decomposition $f = h + \overline{g}$, where $h$ and $g$ are analytic;
this decomposition is unique up to an additive constant (see \cite{CS1984,Duren2004}).
 
It is a classical result that a complex-valued harmonic function (or  harmonic mapping) \( f \) is locally univalent and sense-preserving in \( \Omega \) if and only if its Jacobian determinant
$$
J_f:= |f_z|^2 - |f_{\bar z}|^2>0
$$
in $ \Omega $. Equivalently, the  second complex dilatation  $\omega_f$ of $ f $ satisfies
$$
|\omega_{f}(z)|
=
\left| \frac{f_{\bar z}(z)}{f_z(z)} \right|
< 1,
\qquad z \in \Omega.
$$

For a domain $\Omega \subset \mathbb{C}$, 
let $\mathcal{A}(\Omega)$ denote the set of all analytic functions from $\Omega$ into $\mathbb{C}$, 
and let $\mathcal{H}(\Omega)$ denote the set of all complex-valued harmonic functions from $\Omega$ into $\mathbb{C}$.
For $f \in \mathcal{H}(\Omega)$, if $J_f \geq 0$ in $\Omega$ and there exist constants $K \geq 1$ and $K_0 \geq 0$ such that
\[
    \Lambda_f^2 \le K J_f + K_0 \quad \text{in } \Omega,
\]
then $f$ is called a  {\it weak harmonic $(K, K_0)$-quasiregular  mapping} (cf. \cite{CK,FS-1958,N-1953}). 
We denote by $\mathcal{S}_{h}(K, K_{0})$ the set of all weak harmonic $(K, K_{0})$-quasiregular  mappings defined in 
$\mathbb{D}$. Obviously, $\mathcal{A}(\mathbb{D})\subset\mathcal{S}_{h}(1, 0)$.

For convenience, we make a notational convention.
Throughout this paper,
 we use the symbol $C$ to denote various positive
constants, whose values may change from one occurrence to another.
Also we denote by $C=C(a_{1},a_{2},\ldots)$ a constant that depends only on the given parameters $a_{1}$, $a_{2}$, $\ldots$
and whose value may vary from one occurrence to another.

The {\it  Hardy type space}
$\mathscr{H}^{p}$ $(p\in(0,\infty])$ consists of all those functions
$f$ of $\mathbb{D}$ into $\mathbb{C}$ such that $f$ is measurable, $M_{p}(r,f)$ exists for all $r\in[0,1)$,
$$\|f\|_{p}:=\sup_{r\in[0,1)}M_{p}(r,f)<\infty$$ for
$p\in(0,\infty)$, and $$\|f\|_{\infty}:=\sup_{r\in[0,1)}M_{\infty}(r,f)<\infty$$ for $p=\infty$,
where
$$M_{p}(r,f)=\left(\frac{1}{2\pi}\int_{0}^{2\pi}|f(re^{i\theta})|^{p}d\theta\right)^{\frac{1}{p}}~\mbox{and}~M_{\infty}(r,f)=\sup_{\theta\in[0,2\pi]}|f(re^{i\theta})|.$$
 In particular, we use $\mathbf{h}^{p}:=\mathscr{H}^{p}\cap\mathcal{H}(\mathbb{D})$ and
$H^{p}:=\mathscr{H}^{p}\cap\mathcal{A}(\mathbb{D})$
to denote the {\it harmonic Hardy space} and the {\it analytic Hardy space}, respectively. If
$f\in\mathbf{h}^{p}$ for some $p\in (1,\infty)$,
then the radial limits
\be\label{eq-j-1}f^{\ast}(\zeta)=\lim_{r\rightarrow1^{-}}f(r\zeta)\ee exist for almost every $\zeta\in\mathbb{T}$, and $f^{\ast}\in L^{p}(\mathbb{T})$ (see \cite[Theorems 6.7,  6.13 and 6.39]{ABR}).
  Since $|f|^{p}$ is subharmonic for $p>1$, we see that $M_{p}(r,f)$ is increasing on $r\in[0,1)$, and
 \be\label{fx-1}
 \|f\|_{p}^{p}=\lim_{r\rightarrow1^{-}}M_{p}^{p}(r,f)=\frac{1}{2\pi}\int_{0}^{2\pi}|f(e^{i\theta})|^{p}d\theta.
 \ee

\subsection{The sharp Littlewood-Paley type constant problem}

Let $d\sigma$ denote the normalized arclength measure on $\mathbb T$, namely
$$
d\sigma(e^{i\theta})=\frac{d\theta}{2\pi}.
$$
For $p\in(0,\infty )$, denote by $L^{p}(\mathbb{T})$ the
set of all measurable functions $\psi$ of $\mathbb{T}$ into
$\mathbb{C}$ with
$$\|\psi\|_{L^{p}}=\left(\int_{0}^{2\pi}|\psi(e^{i\theta})|^{p}d\sigma(e^{i\theta})\right)^{\frac{1}{p}}<\infty.$$

Let $F\in L^p(\mathbb T)$, where   $p\in(1,\infty)$. Then its $n$-th Fourier coefficient is defined by
\[
\widehat F(n)
=
\int_{\mathbb T}F(\zeta)\,\overline{\zeta}^{\,n}\,d\sigma(\zeta),
\qquad n\in\mathbb Z.
\]
Equivalently,
\[
\widehat F(n)
=
\frac{1}{2\pi}
\int_0^{2\pi}
F(e^{i\theta})e^{-in\theta}\,d\theta,
\qquad n\in\mathbb Z.
\]
Thus the Fourier series of $F$ is formally written as
\[
F(\zeta)
=
\sum_{n=-\infty}^{\infty}\widehat F(n)\zeta^n,
\qquad \zeta\in\mathbb T.
\]

 For an integer set $I\subset\mathbb{Z}$, we define the Fourier projection associated with $I$ by
\[
M_I F(\zeta)
=
\sum_{n\in I}\widehat F(n)\zeta^n,
\qquad \zeta\in\mathbb T.
\]

For $k\geq1$, put
\[
I_k^+=[2^{k-1},2^k]\cap\mathbb Z,
\qquad
I_k^-=[-2^k,-2^{k-1}]\cap\mathbb Z.
\]
The dyadic Littlewood-Paley square function is defined by
\[
\gamma(F)(\zeta)
=
\left(
\sum_{k=1}^{\infty}
\left(
|M_{I_k^+}F(\zeta)|^2
+
|M_{I_k^-}F(\zeta)|^2
\right)
\right)^{1/2},
\qquad \zeta\in\mathbb T.
\]

A well-known theorem of Littlewood and Paley \cite{LP1936} asserts that
\begin{equation}\label{eq-1.1t}
\|\gamma(F)\|_p \leq A_p \|F\|_p,
\end{equation}
where the positive constants $A_p$ depend only on $p$.

In 1989, Bourgain established the sharp endpoint behavior, which is stated as follows.

\begin{Thm} {\rm (\cite[Theorem 1]{Bourgain1989})} \label{Thm-1}
For $p\in(1,2]$, let $F\in L^p(\mathbb T)$. Then
there exist absolute constants $0<C_{1},~C_{2}<\infty$ such  that
$$C_{1}(p-1)^{-3/2}<A_p<C_{2}(p-1)^{-3/2},$$ where $A_p$ is the constant for which inequality (\ref{eq-1.1t}) holds.
\end{Thm}

\begin{Thm} {\rm (\cite[Theorem 2]{Bourgain1989})} \label{Thm-2}
For $p\in(2,\infty)$, let $F\in L^p(\mathbb T)$. Then
there exist absolute constants $0<C_{3},~C_{4}<\infty$ such  that
$$C_{3}p<A_p<C_{4}p,$$ where $A_p$ is the constant for which inequality (\ref{eq-1.1t}) holds.
\end{Thm}

%In 1989, Bourgain \cite{Bourgain1989} determined the asymptotic behavior of $A_p$ by showing
%\beqq\label{eq-1.1h}
%A_p = O(p),\quad p\to\infty,\qquad \text{and} \qquad A_p = O\bigl((p-1)^{-3/2}\bigr),\quad p\to 1^+.
%\eeqq

By replacing the condition $ F \in L^{p}(\mathbb T) $ with $ f \in H^{p} $, Pichorides~\cite{Pichorides1992} obtained substantially sharper asymptotic estimates for $ p \in (1,2] $, as detailed below.

\begin{Thm} {\rm (\cite[Theorem]{Pichorides1992})} \label{Thm-A}
For $p\in(1,2]$, let $f\in H^p$. Then  there exists a positive constant $A_p$
such that
$$\|\gamma(f^{\ast})\|_p \leq A_p \|f\|_p,$$
where  $f^{\ast}$ is defined in (\ref{eq-j-1}) and $A_p = O\bigl((p-1)^{-1}\bigr)$ as $p \to 1^+$.
\end{Thm}

However, for $ p > 2 $ and $ f \in H^{p} $, Pichorides demonstrated that the sharp growth order of $ A_p $ remains identical to that in Theorem~\ref{Thm-2} as $ p \to \infty $.
For further related work on Littlewood-Paley type inequalities, we refer the reader to \cite{Bakas, Rubio1985, Xu2022}.  
Motivated by the aforementioned studies, we investigate Littlewood-Paley type inequalities for harmonic mappings and 
establish a sharp endpoint behavior, as presented below.

%\begin{Thm}\cite[Theorem]{Pichorides1992}\label{Thm-A}
%If $f \in H^p(\mathbb{D})$, $p>1$, then \eqref{eq-1.1t} holds with $A_p = O\bigl((p-1)^{-1}\bigr)$ as $p \to 1^+$.
%\end{Thm}

\begin{thm}\label{thm-1}
Let $p\in(1,2]$, $K\geq 1$ and $K_{0}\geq 0$. Suppose $f=h+\overline{g}$ is a  weak harmonic $(K,K_{0})$-quasiregular type mapping in 
$\mathbb{D}$, where $h$ and $g$ are analytic with $g(0)=0$. 
If $f\in\mathbf{h}^{p}$, then  there exists  positive constants $C(p,K)$ and $C‘(p,K)$ independent of $f$ such that
$$
\|\gamma(f^{\ast})\|_{L^{p}(\mathbb{T})} \leq C(p,K)\, \|f\|_{p}+\sqrt{K_{0}}C’(p,K),
$$ where $f^{\ast}$ is defined in (\ref{eq-j-1}). Moreover, the constants can be chosen so that
\[
C(p,K) = O\bigl((p-1)^{-1}\bigr), \qquad
C'(p,K) = O\bigl((p-1)^{-1}\bigr)
\]
as $p \to 1^+$. For each fixed $K\geq 1$ and $K_0\geq 0$, the order
$(p-1)^{-1}$ of the coefficient of $\|f\|_p$ is sharp.
\end{thm}

\begin{rem}
Since $\mathcal{A}(\mathbb{D}) \subset \mathcal{S}_{h}(1, 0)$, Theorem~\ref{thm-1} improves and generalizes Theorems~\ref{Thm-A}.
\end{rem}

\subsection{The sharp Hardy-Littlewood type constant problem}

For every fixed $p \in (0,\infty]$ and each integer $n \geq 1$, let $C(n,p)$ denote the smallest constant such that
\begin{equation}
\label{eq:main}
|a_n| + |b_n| \leq C(n,p) \| f \|_p
\end{equation}
for every complex-valued harmonic function
\[
f(z) = h(z) + \overline{g(z)} \in \mathbf{h}^p, 
\]
where 
$
h(z) = \sum_{n=0}^{\infty} a_n z^n$ and $ g(z) = \sum_{n=1}^{\infty} b_n z^n.
$

When $g \equiv 0$ and $p \in [1,\infty)$, the triangle inequality combined with Hölder's inequality readily yields $C(n,p) = 1$ for all $n \geq 1$. For $g \equiv 0$ and $p = \infty$, Cauchy's coefficient estimate gives $C(n,\infty) = 1$ for all $n \geq 1$. In both cases, the constant is sharp, as demonstrated by the extremal function $f(z) = z^n$.

For $g \equiv 0$ and $p \in (0,1)$, Hardy and Littlewood \cite{HardyLittlewood} first established the general estimate
\begin{equation}
C(n,p) \leq C_p \, n^{1/p - 1}, \qquad n \geq 1,
\end{equation}
where $C_p \geq 1$ depends only on $p$. The first exact value was later obtained by Bondarenko, Brevig, Saksman, and Seip \cite{BondarenkoBrevigSaksmanSeip2019}, who proved that
\begin{equation}
C(1,p) = \sqrt{\frac{2}{p}} \left(1 - \frac{p}{2}\right)^{\frac{1}{p} - \frac{1}{2}},
\end{equation}
where $g \equiv 0$ and $p \in (0,1)$.

Subsequently, Brevig and Saksman \cite{BrevigSaksman2020} determined
\begin{equation}
C(2,p) = \frac{2}{p} \left(1 - \frac{p}{2}\right)^{\frac{2}{p} - 1},
\end{equation}
and also computed the special value
\begin{equation}
C\left(3,\frac{2}{3}\right) = \sqrt{\frac{2\left(1103 + 33\sqrt{33}\right)}{1153}}, 
\end{equation}
where again $g \equiv 0$ and $p \in (0,1)$.

For all other positive integers $n$, the exact value of $C(n,p)$ remains an open problem, still under the assumption $g \equiv 0$ and $p \in (0,1)$.

However, the optimal constant $C(n,p)$ exhibits a significant difference between the cases $g \equiv 0$ and $g \not\equiv 0$. 
For example, in the case $g \not\equiv 0$, it is known that $C(n,\infty)=4/\pi$. 
When $g \not\equiv 0$ and $p\in(0,\infty)$, 
the exact value of $C(n,p)$ remains an open problem. 
In particular, Chen, Ponnusamy and Wang \cite{cpw-2012} established an upper bound for $C(n,p)$ in the setting $g \not\equiv 0$ and $p\in[1,\infty)$. 
Their result is stated as follows.

\begin{Thm}{\rm (\cite[Theorem 3]{cpw-2012})}\label{Thm-A1}
Let $f$ be a harmonic mapping in $\mathbb{D}$ such that
\beqq
f(z)
=
\sum_{n=0}^{\infty}a_nz^n
+
\sum_{n=1}^{\infty}\overline{b_n}\,\overline{z}^{\,n},
\eeqq
and $f\in\mathbf{h}^p$ for some $p\in[1,\infty]$. Then the
following statements hold:
\begin{enumerate}
    \item
    $
    |a_0|\leq \|f\|_p.
    $
    \item For $p\in[1,\infty)$,
    \beqq
    |a_n|+|b_n|
    \leq
    \frac{
    2^{\frac1p+2}(1+np)^{n+\frac1p}
    }{
    \pi(pn)^n
    }
    \|f\|_p,
    \qquad n\geq1.
    \eeqq
    \item For $p=\infty$,
    \beqq
    |a_n|+|b_n|
    \leq
    \frac4\pi\|f\|_\infty,
    \qquad n\geq1.
    \eeqq
    The estimate in this case is sharp, and the  extremal
    functions are
    \beqq
    f_n(z)
    =
    \frac{2\alpha}{\pi}\|f\|_\infty
    \arg\left(
    \frac{1+\beta z^n}{1-\beta z^n}
    \right),
    \eeqq
    where
   $
    |\alpha|=|\beta|=1
    $ and ``$\arg$'' denotes the principal value of the argument.
\end{enumerate}
\end{Thm}

The second purpose of this paper is to determine the exact value of $C(n,p)$ for $p\in[1,\infty)$. Our  result  is stated as follows.

\begin{thm}\label{thm-2}
Let $f$ be a harmonic mapping in $\mathbb{D}$ such that
\beqq
f(z)
=
\sum_{n=0}^{\infty}a_nz^n
+
\sum_{n=1}^{\infty}\overline{b_n}\,\overline{z}^{\,n},
\eeqq
and $f\in\mathbf{h}^p$ for some $p\in[1,\infty)$. Then,  for every $n\geq1$,
\beqq
|a_n|+|b_n|
\leq C(n,p)\|f\|_p,
\eeqq
where $
C(n,1)=2
$
and
\beqq
C(n,p)
=
2
\left[
\frac{
\Gamma\left(\frac{2p-1}{2(p-1)}\right)
}{
\sqrt{\pi}\,
\Gamma\left(\frac{3p-2}{2(p-1)}\right)
}
\right]^{\frac{p-1}{p}}
\eeqq
for $1<p<\infty$.
All these estimates are sharp. 
\end{thm}

The proofs of Theorems \ref{thm-1} and \ref{thm-2} will be presented in Section \ref{sec-2}.

\section{The proofs of the main results}\label{sec-2}

Prior to proving Theorem \ref{thm-1}, we shall recall some necessary Lemmas.

%The following result follows from \cite[Theorem  1.1]{Pav-2013}.

 %It follows from (\ref{eq-1.1}) and (\ref{eq-1.2}) that, for $z\in\mathbb{D}$, 

%\beqq
%|F_{2}'|\leq\Lambda_{f}\leq\frac{2K}{1+K}|h'|+\frac{\sqrt{K_{0}}}{1+K}
%\eeqq and 

%\beqq
%|F_{1}'|\geq\lambda_{f}\leq\frac{2}{1+K}|h'|-\frac{\sqrt{K_{0}}}{1+K},
%\eeqq
%which imply that

%\beqq
%|F_{2}'|\leq\,K|F_{1}'|+\sqrt{K_{0}}.
%\eeqq

 %\begin{Lem}\label{Lemx}
%Suppose that $a,~b\in[0,\infty)$ and $q\in(0,\infty)$. Then
%$$(a+b)^{q}\leq2^{\max\{q-1,0\}}(a^{q}+b^{q}).$$
%\end{Lem}

\begin{Lem}\label{Lemx} {\rm (\cite[p. 158]{Kuang2010})} For $n \ge 2$, let $a_1, \dots, a_n \in \mathbb{C}$. Then
\[
\left( \sum_{k=1}^n |a_k| \right)^p \le
\begin{cases}
	\sum_{k=1}^n |a_k|^p, & \text{if } 0 < p < 1, \\[4pt]
	n^{p-1} \sum_{k=1}^n |a_k|^p, & \text{if } p \ge 1.
\end{cases}
\]
\end{Lem}

\begin{Lem}{\rm (\cite[p.1]{Pav2009})}\label{Lem-A}
	If $g \in C^2(\mathbb{D})$, then, for $r\in(0,1)$,
	$$
	\frac{1}{2\pi}\int_{0}^{2\pi} g(re^{i\theta}) d\theta = g(0) + \frac{1}{2}\int_{|z|<r}
 \Delta \left(g(z)\right) \log\frac{r}{|z|} dA(z),
	$$ where $dA(z)=dxdy/\pi$.
\end{Lem}

\subsection*{Proof of Theorem \ref{thm-1}}
	Since $f$ is a  weak harmonic $(K, K_{0})$-quasiregular mapping in $\mathbb{D}$, we have
\begin{equation}\label{eq-3.1}
    \Lambda_f^{2} \leq K\Lambda_f \lambda_f + K_{0}.
\end{equation}

We first prove that
\begin{equation}\label{eq-3.2}
    \Lambda_f^{2}
    \leq
    \frac{2K^{2}}{K^{2}+1}
    \bigl( |h'|^{2} + |g'|^{2} \bigr)
    +
    \frac{2K_{0}}{K^{2}+1}.
\end{equation}

Indeed, since
\[
|h'|^{2} + |g'|^{2}
=
\frac{\Lambda_f^{2} + \lambda_f^{2}}{2},
\]
we obtain
\begin{align*}
    \Lambda_f^{2}
    - \frac{2K^{2}}{K^{2}+1}
    \bigl( |h'|^{2} + |g'|^{2} \bigr)
    &=
    \Lambda_f^{2}
    - \frac{K^{2}}{K^{2}+1}
    \bigl( \Lambda_f^{2} + \lambda_f^{2} \bigr) \\
    &=
    \frac{
        \bigl( \Lambda_f - K\lambda_f \bigr)
        \bigl( \Lambda_f + K\lambda_f \bigr)
    }{K^{2}+1}.
\end{align*}

If $\Lambda_f - K\lambda_f \leq 0,$ then \eqref{eq-3.2} follows immediately.  
If $\Lambda_f - K\lambda_f > 0$, then $K\lambda_f < \Lambda_f$, and therefore
\begin{align*}
    \bigl( \Lambda_f - K\lambda_f \bigr)
    \bigl( \Lambda_f + K\lambda_f \bigr)
    &\leq
    2\Lambda_f \bigl( \Lambda_f - K\lambda_f \bigr) \\
    &=
    2\bigl( \Lambda_f^{2} - K\Lambda_f \lambda_f \bigr) \\
    &\leq 2K_{0},
\end{align*}
where the last inequality follows from \eqref{eq-3.1}.  
Thus \eqref{eq-3.2} holds in both cases.
	
	For $\varepsilon>0$, define
	\[
	V_{\varepsilon}(z)
	=
	\left(
	|f(z)|^{2}
	+\frac{4K_{0}}{K^{2}+1}|z|^{2}
	+\varepsilon
	\right)^{1/2}
	\]
	and
	\[
	U_{\varepsilon}(z)
	=
	\left(
	|h(z)|^{2}
	+|g(z)|^{2}
	+\frac{4K_{0}}{K^{2}+1}|z|^{2}
	+\varepsilon
	\right)^{1/2}.
	\]
	We next compare $\Delta(U_{\varepsilon}^{p})$ with
	$\Delta(V_{\varepsilon}^{p})$.
	
	Let
\[
\Phi(z)= \bigl(\Phi_{1}(z), \Phi_{2}(z)\bigr) \in \mathbb{C}^{2} \cong \mathbb{R}^{4},
\]
where $\Phi_{1}(z)=f(z)$ and $\Phi_{2}(z)=2\sqrt{\dfrac{K_{0}}{K^{2}+1}}\,z$. Then
\[
|\Phi(z)|^{2}=|f(z)|^{2}+\frac{4K_{0}}{K^{2}+1}|z|^{2},
\]
and
\begin{equation}\label{eq-3.3}
\Delta\left(|\Phi(z)|^{2}\right)
=4\left( |h'|^{2}+|g'|^{2}+\frac{4K_{0}}{K^{2}+1} \right).
\end{equation}

For $0\leq \vartheta \leq 2\pi$, set
\[
\xi_{\vartheta}=(\cos\vartheta,\sin\vartheta)\in \mathbb{R}^{2}.
\]
Under the identification $\mathbb{R}^{2}\cong \mathbb{C}$, the vector $\xi_{\vartheta}$ corresponds to $e^{i\vartheta}$. Hence,
\begin{align*}
D\Phi(z)\xi_{\vartheta}
&=\lim_{t\to 0} \frac{\Phi(z+t e^{i\vartheta})-\Phi(z)}{t} \\
&=
\left(
\partial_{\vartheta} f(z),\,
2\sqrt{\frac{K_{0}}{K^{2}+1}}\, e^{i\vartheta}
\right),
\end{align*}
where
\[
D\Phi(z)
=
\begin{pmatrix}
\dfrac{\partial \Phi_{1}}{\partial x}(z) & \dfrac{\partial \Phi_{1}}{\partial y}(z) \\[1.2ex]
\dfrac{\partial \Phi_{2}}{\partial x}(z) & \dfrac{\partial \Phi_{2}}{\partial y}(z)
\end{pmatrix}.
\]
Here $D\Phi:~\mathbb{R}^{2}\rightarrow\mathbb{R}^{4}$ acts on real vectors.
Consequently,
\[
|D\Phi(z)\xi_{\vartheta}|^{2}
=
|\partial_{\vartheta} f(z)|^{2}
+
\frac{4K_{0}}{K^{2}+1}.
\]
	By the definition of the operator norm, we have
\begin{align*}
  \|D\Phi(z)\|^{2}
  &=
  \max_{|\xi|=1} |D\Phi(z)\xi|^{2} \\
  &=
  \max_{0\leq\vartheta\leq2\pi}
  \left(
    \left|\partial_{\vartheta}f(z)\right|^{2}
    +
    \frac{4K_{0}}{K^{2}+1}
  \right) \notag \\
  &=
  \left(
    \max_{0\leq\vartheta\leq2\pi}
    \left|\partial_{\vartheta}f(z)\right|
  \right)^{2}
  +
  \frac{4K_{0}}{K^{2}+1} \notag \\
  &=
  (\Lambda_f(z))^{2}
  +
  \frac{4K_{0}}{K^{2}+1}.
\end{align*}

Comparing these identities, we obtain
\begin{equation*}
  \nabla\left(|\Phi(z)|^{2}\right)
  =
  2\Bigg(
    \operatorname{Re}\!\left(
      \frac{\partial\phi_{1}}{\partial x}\overline{\phi_{1}}
      +
      \frac{\partial\phi_{2}}{\partial x}\overline{\phi_{2}}
    \right),
    \operatorname{Re}\!\left(
      \frac{\partial\phi_{1}}{\partial y}\overline{\phi_{1}}
      +
      \frac{\partial\phi_{2}}{\partial y}\overline{\phi_{2}}
    \right)
  \Bigg)
  =
  2\operatorname{Re}\!\left((D\Phi(z))^{T}\overline{\Phi(z)}\right),
\end{equation*}
and
\begin{equation*}
  V_{\varepsilon}^{2}(z)
  =
  |\Phi(z)|^{2}+\varepsilon
  \geq
  |\Phi(z)|^{2},
\end{equation*}
which together imply that
\begin{equation}\label{eq-3.4}
  \left|\nabla\left(|\Phi(z)|^{2}\right)\right|^{2}
  \leq
  4\|D\Phi(z)\|^{2}\,|\Phi(z)|^{2}
  \leq
  4V_{\varepsilon}^{2}(z)
  \left(
    (\Lambda_f(z))^{2}
    +
    \frac{4K_{0}}{K^{2}+1}
  \right),
\end{equation}
where the superscript ``$T$'' denotes the transpose of a matrix.
Here ``$\nabla$'' denotes the real gradient. 

	Since \(1 < p \leq 2\), applying the chain rule along with \eqref{eq-3.3} and \eqref{eq-3.4} yields
\beq\label{eq-3.5}
	\Delta(V_{\varepsilon}^{p}) &=& \frac{p}{2}V_{\varepsilon}^{p-2}\Delta(V_{\varepsilon}^{2}) +
	\frac{p}{2}\left(\frac{p}{2}-1\right)V_{\varepsilon}^{p-4}|\nabla V_{\varepsilon}^{2}|^{2} \\ \nonumber
	&\geq&
	2pV_{\varepsilon}^{p-2}
	\left(
	|h'|^{2}
	+|g'|^{2}
	+\frac{4K_{0}}{K^{2}+1}
	\right)
	+
	p(p-2)V_{\varepsilon}^{p-2}
	\left(
	\Lambda_f^{2}
	+\frac{4K_{0}}{K^{2}+1}
	\right)\\ \nonumber
	&=&
	pV_{\varepsilon}^{p-2}
	\psi(g,h),\eeq
where $$\psi(g,h)=2\left(|h'|^{2}+|g'|^{2}\right)
	-(2-p)\Lambda_f^{2}
	+\frac{4pK_{0}}{K^{2}+1}.$$
We now turn to estimating $\psi(g,h)$.
From \eqref{eq-3.2}, it follows that
\begin{equation}\label{eq-3.6}
	\psi(g,h)
	\geq
	\left[
	2-\frac{2(2-p)K^{2}}{K^{2}+1}
	\right]
	\left(
	|h'|^{2}+|g'|^{2}
	\right)+
	\frac{(6p-4)K_{0}}{K^{2}+1}.
\end{equation}

Since, for $1 < p \leq 2$,  we have
\[
2-\frac{2(2-p)K^{2}}{K^{2}+1}
\geq
\frac{2}{K^{2}+1}
\]
and
\[
\frac{(6p-4)K_{0}}{K^{2}+1}
\geq
\frac{2K_{0}}{K^{2}+1}
=
\frac{1}{2}\cdot
\frac{4K_{0}}{K^{2}+1},
\]
it follows from \eqref{eq-3.5} and \eqref{eq-3.6} that
\begin{equation}\label{eq-3.7}
	\Delta(V_{\varepsilon}^{p})
	\geq
	p\min\left\{
	\frac{2}{K^{2}+1},
	\frac{1}{2}
	\right\}
	V_{\varepsilon}^{p-2}
	\left(
	|h'|^{2}
	+|g'|^{2}
	+\frac{4K_{0}}{K^{2}+1}
	\right).
\end{equation}
	
	On the other hand, direct differentiation gives
	\begin{align*}
		\Delta(U_{\varepsilon}^{p})
		={}&
		2pU_{\varepsilon}^{p-2}
		\left(
		|h'|^{2}
		+|g'|^{2}
		+\frac{4K_{0}}{K^{2}+1}
		\right)\\
		&+
		p(p-2)U_{\varepsilon}^{p-4}
		\left|
		h'\overline{h}
		+g'\overline{g}
		+\frac{4K_{0}}{K^{2}+1}\overline{z}
		\right|^{2},
	\end{align*}
	which, together with $p-2\leq0$, yields that
	\begin{equation}\label{eq-3.8}
		\Delta(U_{\varepsilon}^{p})
		\leq
		2pU_{\varepsilon}^{p-2}
		\left(
		|h'|^{2}
		+|g'|^{2}
		+\frac{4K_{0}}{K^{2}+1}
		\right).
	\end{equation}
	Furthermore,
\[
|f|^{2}
=
|h+\overline{g}|^{2}
\leq
2\bigl(|h|^{2}+|g|^{2}\bigr),
\]
and hence
\[
V_{\varepsilon}^{2}\leq 2U_{\varepsilon}^{2}.
\]
Since \(p-2\leq 0\), this inequality yields
\[
U_{\varepsilon}^{p-2}
\leq
2^{(2-p)/2}\, V_{\varepsilon}^{p-2}.
\]
	Combining  \eqref{eq-3.7} and
	\eqref{eq-3.8} gives
	\begin{equation}\label{eq-3.9}
		\Delta(U_{\varepsilon}^{p})
		\leq
		\frac{2^{1+(2-p)/2}}
		{\displaystyle
			\min\left\{
			\frac{2}{K^{2}+1},\frac{1}{2}
			\right\}}
		\Delta(V_{\varepsilon}^{p})\leq
		\frac{2^{3/2}}
		{\displaystyle
			\min\left\{
			\frac{2}{K^{2}+1},\frac{1}{2}
			\right\}}
		\Delta(V_{\varepsilon}^{p}).
	\end{equation}
	
	From the assumption $g(0)=0$, we have $f(0)=h(0)$ and therefore
	\[
	U_{\varepsilon}(0)=V_{\varepsilon}(0).
	\]
	Applying Lemma~\ref{Lem-A} to \eqref{eq-3.9}, we obtain,
	for every $0<r<1$,
	$$
		\int_{\mathbb{T}}
		U_{\varepsilon}(r\zeta)^{p}
		\,d\sigma(\zeta)
		-U_{\varepsilon}(0)^{p}\leq
		\frac{2^{3/2}}
		{\displaystyle
			\min\left\{
			\frac{2}{K^{2}+1},\frac{1}{2}
			\right\}}
		\left[
		\int_{\mathbb{T}}
		V_{\varepsilon}(r\zeta)^{p}
		\,d\sigma(\zeta)
		-V_{\varepsilon}(0)^{p}
		\right].
	$$
	Since
	$$
	\frac{2^{3/2}}
		{\min\left\{
		\frac{2}{K^{2}+1},\frac{1}{2}
		\right\}}
	\geq1,
	$$
	the equality $U_{\varepsilon}(0)=V_{\varepsilon}(0)$ implies that
	\begin{align*}
		\int_{\mathbb{T}}
		U_{\varepsilon}(r\zeta)^{p}
		\,d\sigma(\zeta)
		\leq
		\frac{2^{3/2}}
		{\displaystyle
			\min\left\{
			\frac{2}{K^{2}+1},\frac{1}{2}
			\right\}}
		\int_{\mathbb{T}}
		V_{\varepsilon}(r\zeta)^{p}
		\,d\sigma(\zeta).
	\end{align*}
	Letting $\varepsilon\to0^{+}$ gives
	\be\label{eq-3.10}
		\int_{\mathbb{T}}
		\left(
		|h(r\zeta)|^{2}
		+|g(r\zeta)|^{2}
		+\frac{4K_{0}r^{2}}{K^{2}+1}
		\right)^{p/2}
		\,d\sigma(\zeta)
		\leq
		C(K)
		\int_{\mathbb{T}}
		\left(
		|f(r\zeta)|^{2}
		+\frac{4K_{0}r^{2}}{K^{2}+1}
		\right)^{p/2}
		\,d\sigma(\zeta),
	\ee where $$C(K)=\frac{2^{3/2}}
		{\displaystyle
			\min\left\{
			\frac{2}{K^{2}+1},\frac{1}{2}
			\right\}}.$$
	Since $p/2\leq1$, by Lemma \ref{Lemx}  and \eqref{eq-3.10}, we see that
	\begin{align}\label{eq-3.11}
		\int_{\mathbb{T}}
		\left(
		|h(r\zeta)|^{2}
		+|g(r\zeta)|^{2}
		\right)^{p/2}
		\,d\sigma(\zeta)\leq
		C(K)
		\left[
		\int_{\mathbb{T}}
		|f(r\zeta)|^{p}
		\,d\sigma(\zeta)
		+
		\left(
		\frac{4K_{0}}{K^{2}+1}
		\right)^{p/2}
		\right].
	\end{align}
	Taking the supremum over $r\in(0,1)$ in
	\eqref{eq-3.11}, we conclude that $h,g\in H^{p}$.
	Passing to radial boundary values and then taking the $p$-th root
	gives
	\begin{equation}\label{eq-3.12}
		\left\|
		\left(
		|h^{\ast}|^{2}+|g^{\ast}|^{2}
		\right)^{1/2}
		\right\|_{L^{p}(\mathbb{T})}
		\leq
		\left(C(K)\right)^{\frac{1}{p}}
		\left(
		\|f^{\ast}\|_{L^{p}(\mathbb{T})}
		+C(K_{0})
		\right),
	\end{equation}
where $$C(K_{0})=
		\frac{2\sqrt{K_{0}}}{\sqrt{K^{2}+1}},$$
$$h^{\ast}(\zeta)=\lim_{r\rightarrow1^{-}}h(r\zeta)$$ and $$g^{\ast}(\zeta)=\lim_{r\rightarrow1^{-}}g(r\zeta).$$
	
	Let $S_+$ denote the analytic Littlewood-Paley square function. Then, by the
	definition of $\gamma(f^*)$, we have
	$$
	\left(\gamma(f^*)(\zeta)\right)^2
	=
	\left(S_+(h^*)(\zeta)\right)^2+\left(S_+(g^*)(\zeta)\right)^2,
	\qquad \zeta\in\mathbb T.
	$$ %where $$h^{\ast}(\zeta)=\lim_{r\rightarrow1^{-}}h(r\zeta)$$ and $$g^{\ast}(\zeta)=\lim_{r\rightarrow1^{-}}g(r\zeta).$$
	By Lemma \ref{Lemx}, we obtain
	\beqq
	\|\gamma(f^*)\|_{L^p(\mathbb T)}^p
	&=&\int_{\mathbb T}
	\left[
	\left(S_+(h^*)(\zeta)\right)^2+\left(S_+(g^*)(\zeta)\right)^2
	\right]^{p/2}
	\,d\sigma(\zeta)\\
	&\leq&\int_{\mathbb T}\left(S_+(h^*)(\zeta)\right)^p\,d\sigma(\zeta)
	+\int_{\mathbb T}\left(S_+(g^*)(\zeta)\right)^p\,d\sigma(\zeta)\\
	&=&\|S_+(h^*)\|_{L^p(\mathbb T)}^p
	+\|S_+(g^*)\|_{L^p(\mathbb T)}^p,
	\eeqq 
	which, together with   Theorem \ref{Thm-A} and Lemma \ref{Lemx}, implies that there exist a positive constant $A_{p}$ such that  
	\beqq
	\|\gamma(f^*)\|_{L^p(\mathbb T)}
	\leq 
	A_{p}
	\left(\|h^*\|_{L^p(\mathbb T)}^p
	+\|g^*\|_{L^p(\mathbb T)}^p\right)^{\frac{1}{p}}\leq A_{p}2^{1/p-1/2}
	\left\|
	\left(
	|h^{\ast}|^{2}
	+
	|g^{\ast}|^{2}
	\right)^{1/2}
	\right\|_{L^{p}(\mathbb{T})} .
	\eeqq
	Combining this estimate with \eqref{eq-3.12}, we obtain
	\beqq
	\|\gamma(f^{\ast})\|_{L^{p}(\mathbb{T})}
	&\leq&
	\left(C(K)\right)^{\frac{1}{p}}2^{1/p-1/2}A_{p}
	\left(
	\|f^{\ast}\|_{L^{p}(\mathbb{T})}
	+C(K_{0})
	\right)\\
&=&\left(C(K)\right)^{\frac{1}{p}}2^{1/p-1/2}A_{p}
	\left(
	\|f\|_{p}
	+C(K_{0})
	\right).
	\eeqq
	It remains to prove that the order of $C(p,K)$ is sharp. Let $F\in H^p$.
	Since $tF\in\mathcal{S}_{h}(K,K_0)$ for every $t>0$, the estimate in the theorem gives
	$$
	t\|\gamma(F^\ast)\|_{L^p(\mathbb T)}
	\leq t C(p,K)\|F\|_p+\sqrt{K_{0}}C'(p,K).
	$$
	Dividing by $t$ and letting $t\to\infty$, we obtain
	$$
	\|\gamma(F^\ast)\|_{L^p(\mathbb T)}\leq C(p,K)\|F\|_p.
	$$

	Following the construction in the proof of \cite[Proposition 6]{BakasRodriguezLopezSola2019},
	we use a modulated de la Vall\'ee Poussin kernel. For each integer $m\geq1$, let
	$$
	Q_m(\zeta)=\frac1m\left|\sum_{j=0}^{m-1}\zeta^j\right|^2,
	\qquad \zeta\in\mathbb T.
	$$
	The Fourier coefficients of $Q_m$ are $\widehat Q_m(j)=\max\{1-|j|/m,0\}$.
	In particular, $\|Q_m\|_{L^1(\mathbb T)}=1$ and
	$\|Q_m\|_{L^\infty(\mathbb T)}=m$. Hence
	$\|Q_m\|_{L^p(\mathbb T)}\leq m^{(p-1)/p}$.

	Let $N=2^k$ for an integer $k\geq1$. All frequencies of
	$\zeta^{2N}(2Q_{2N}(\zeta)-Q_N(\zeta))$ are positive, so this is the boundary
	function of an analytic polynomial $F_N$. Since the Fourier coefficients of
	$2Q_{2N}-Q_N$ equal $1$ for $|j|\leq N$, we have
	$$
	M_{I_{k+1}^+}F_N^\ast(\zeta)=\sum_{n=N}^{2N}\zeta^n.
	$$
	Since $|\zeta^{2N}|=1$ on $\mathbb T$, it follows from \eqref{fx-1} and
	the triangle inequality that
	$$
	\|F_N\|_p
	=\|F_N^\ast\|_{L^p(\mathbb T)}
	\leq 2(2N)^{(p-1)/p}+N^{(p-1)/p}
	\leq 3(2N)^{(p-1)/p}.
	$$
	To estimate this frequency block, we write $\zeta=e^{it}$. For $0<t<2\pi$,
	the geometric sum formula gives
	$$
	\sum_{n=N}^{2N}e^{int}
	=e^{iNt}\sum_{j=0}^{N}e^{ijt}
	=e^{iNt}\frac{1-e^{i(N+1)t}}{1-e^{it}}.
	$$
	Since $|e^{iNt}|=1$ and
	$1-e^{ix}=-2i e^{ix/2}\sin(x/2)$ for $x\in\mathbb R$, we obtain
	$$
	\left|\sum_{n=N}^{2N}e^{int}\right|
	=\frac{|1-e^{i(N+1)t}|}{|1-e^{it}|}
	=\frac{|\sin((N+1)t/2)|}{|\sin(t/2)|}.
	$$
	By H\"older's inequality and the symmetry of the integrand under
	$t\mapsto2\pi-t$, we obtain
	\begin{align*}
	\|\gamma(F_N^\ast)\|_{L^p(\mathbb T)}
	&\geq \|M_{I_{k+1}^+}F_N^\ast\|_{L^1(\mathbb T)}\\
	&=\frac1\pi\int_0^\pi
	\frac{|\sin((N+1)t/2)|}{\sin(t/2)}\,dt\\
	&\geq \frac2\pi\int_0^{N\pi/2}\frac{|\sin s|}{s}\,ds\\
	&\geq \frac4{\pi^2}\sum_{j=1}^{N/2}\frac1j
	\geq \frac2{\pi^2}\log N.
	\end{align*}
	The second inequality follows from $\sin(t/2)\leq t/2$ and the
	substitution $s=(N+1)t/2$. The third inequality follows by integrating over
	$[(j-1)\pi,j\pi]$, where $s\leq j\pi$ and the integral of $|\sin s|$ equals $2$.
	The last inequality follows from
	$\sum_{j=1}^{N/2}j^{-1}\geq \log(1+N/2)\geq \frac12\log N$.

	Choose $k=\lceil 1/(p-1)\rceil$. Then
	$$
	\log N\geq\frac{\log 2}{p-1}
	\qquad\text{and}\qquad
	(2N)^{(p-1)/p}\leq 2^{(2p-1)/p}\leq 4.
	$$
	Consequently, $\|F_N\|_p\leq12$, and applying the preceding estimate for
	analytic functions to $F_N$ yields
	$$
	C(p,K)\geq
	\frac{\|\gamma(F_N^\ast)\|_{L^p(\mathbb T)}}{\|F_N\|_p}
	\geq\frac{\log 2}{6\pi^2(p-1)}.
	$$
	This completes the proof.
\qed

\begin{Thm}{\rm (\cite[Theorem 1.10]{pavlovic-2019})}\label{ThmB}
For $p\in(1,\infty)$, let  $g\in L^p(\mathbb{T})$.
Suppose that $$f(z)=\frac{1}{2\pi}\int_{0}^{2\pi}\frac{1-|z|^{2}}{|1-ze^{-it}|^{2}}g(e^{it})dt.$$
Then
\beqq
||f||_p=||g||_{L^p(\mathbb{T})}.
\eeqq
\end{Thm}

\subsection*{Proof of Theorem \ref{thm-2}}

Fix an integer $n\geq1$. For $r\in(0,1)$, we have
\beqq
f(re^{it})
=
\sum_{j=0}^{\infty}a_jr^je^{ijt}
+
\sum_{j=1}^{\infty}\overline{b_j}\,r^je^{-ijt},
\eeqq
which yields that
\be\label{cjl-1}
a_nr^n
=
\frac1{2\pi}
\int_0^{2\pi}
f(re^{it})e^{-int}\,dt
\ee
and
\be\label{cjl-2}
\overline{b_n}\,r^n
=
\frac1{2\pi}
\int_0^{2\pi}
f(re^{it})e^{int}\,dt.
\ee
Choose $\lambda,\mu\in\mathbb{T}$ such that
\be\label{jkl-1}
\lambda a_n=|a_n|
\qquad\text{and}\qquad
\mu\overline{b_n}=|b_n|.
\ee
Combining (\ref{cjl-1}),  (\ref{cjl-2}) and (\ref{jkl-1}) gives
\be\label{jkl-2}
r^n\bigl(|a_n|+|b_n|\bigr)
=
\left|
\frac1{2\pi}
\int_0^{2\pi}
f(re^{it})
\bigl(
\lambda e^{-int}+\mu e^{int}
\bigr)\,dt
\right|.
\ee

\noindent $\mathbf{Case~1.}$ Let $p=1$.

 Since
\beqq
\left|
\lambda e^{-int}+\mu e^{int}
\right|
\leq2,
\eeqq
by (\ref{jkl-2}),
we see that
\beqq
\begin{aligned}
r^n\bigl(|a_n|+|b_n|\bigr)
&\leq
\frac1{2\pi}
\int_0^{2\pi}
|f(re^{it})|
\left|
\lambda e^{-int}+\mu e^{int}
\right|\,dt\\
&\leq
\frac1{\pi}
\int_0^{2\pi}
|f(re^{it})|\,dt\\
&=
2M_1(r,f)\\
&\leq
2\|f\|_1.
\end{aligned}
\eeqq
By letting $r\to1^-$, we obtain
\beqq
|a_n|+|b_n|
\leq
2\|f\|_1.
\eeqq

\noindent $\mathbf{Case~2.}$ Let $p\in(1,\infty)$.

 By H\"older's inequality, we have
\beq\label{1.1}
\begin{aligned}
r^n\bigl(|a_n|+|b_n|\bigr)
&\leq
M_p(r,f)
\left(
\frac1{2\pi}
\int_0^{2\pi}
\left|
\lambda e^{-int}+\mu e^{int}
\right|^{\frac{p}{p-1}}\,dt
\right)^{\frac{p-1}{p}}.
\end{aligned}
\eeq
Since $|\lambda|=|\mu|=1$, there exists
$\theta_0\in[0,\pi)$ such that
\beqq
\mu\overline{\lambda}=e^{2i\theta_0}.
\eeqq
Consequently,
\beqq
\begin{aligned}
\left|
\lambda e^{-int}+\mu e^{int}
\right|
&=
\left|
\lambda e^{-int}
\right|
\left|
1+\mu\overline{\lambda}e^{2int}
\right|\\
&=
\left|
1+e^{2i(nt+\theta_0)}
\right|\\
&=
\left|
e^{-i(nt+\theta_0)}
+
e^{i(nt+\theta_0)}
\right|\\
&=
2\left|\cos(nt+\theta_0)\right|,
\end{aligned}
\eeqq which, together with (\ref{1.1}), implies that

\beq\label{1.2}
\left(
\frac1{2\pi}
\int_0^{2\pi}
\left|
\lambda e^{-int}+\mu e^{int}
\right|^{\frac{p}{p-1}}\,dt
\right)^{\frac{p-1}{p}}
=
2\left(
\frac1{2\pi}
\int_0^{2\pi}
\left|\cos(nt+\theta_0)\right|^{\frac{p}{p-1}}\,dt
\right)^{\frac{p-1}{p}}.
\eeq

We now estimate the integral of the right-hand side. Since $ n $ is a positive integer, the change of variables $ u = n t + \theta_0 $, combined with the periodicity of the function 
$ u \mapsto |\cos u|^{\frac{p}{p-1}} $, yields
\beq\label{cjl-3}
\frac1{2\pi}
\int_0^{2\pi}
|\cos(nt+\theta_0)|^{\frac{p}{p-1}}\,dt
&=&
\frac1{2\pi n}
\int_{\theta_0}^{2\pi n+\theta_0}
|\cos u|^{\frac{p}{p-1}}\,du\\ \nonumber
&=&
\frac1{2\pi}
\int_0^{2\pi}
|\cos u|^{\frac{p}{p-1}}\,du\\ \nonumber
&=&
\frac2\pi
\int_0^{\pi/2}
\cos^{\frac{p}{p-1}}u\,du.
\eeq
To evaluate the last integral, we set $s = \sin^2 u$. Then
\beqq
ds=2\sin u\cos u\,du,
\eeqq
which implies that, for $0\leq u\leq\pi/2$,
\beqq
du
=
\frac{ds}{2\sqrt{s}\sqrt{1-s}}
\eeqq
and
\beqq
\cos^{\frac{p}{p-1}}u
=
(1-s)^{\frac{p}{2(p-1)}}.
\eeqq
Consequently,
\beqq
\begin{aligned}
\frac{2}{\pi}
\int_0^{\pi/2}
\cos^{\frac{p}{p-1}}u\,du
&=
\frac{1}{\pi}
\int_0^1
s^{-1/2}
(1-s)^{\frac{1}{2(p-1)}}\,ds\\
&=
\frac{1}{\pi}
B\left(
\frac12,
\frac{2p-1}{2(p-1)}
\right)\\
&=
\frac{1}{\pi}
\frac{
\Gamma\left(\frac12\right)
\Gamma\left(\frac{2p-1}{2(p-1)}\right)
}{
\Gamma\left(\frac{3p-2}{2(p-1)}\right)
}\\
&=
\frac{
\Gamma\left(\frac{2p-1}{2(p-1)}\right)
}{
\sqrt{\pi}\,
\Gamma\left(\frac{3p-2}{2(p-1)}\right)
},
\end{aligned}
\eeqq
which, together with (\ref{cjl-3}), yields that
\beq\label{1.3}
\frac1{2\pi}
\int_0^{2\pi}
|\cos(nt+\theta_0)|^{\frac{p}{p-1}}\,dt
=
\frac1{2\pi}
\int_0^{2\pi}
|\cos u|^{\frac{p}{p-1}}\,du
=
\frac{
\Gamma\left(\frac{2p-1}{2(p-1)}\right)
}{
\sqrt{\pi}\,
\Gamma\left(\frac{3p-2}{2(p-1)}\right)
}.
\eeq
Here $B(\cdot,\cdot)$ and $\Gamma(\cdot)$ denote the Beta function and the Gamma function, respectively.
Combining \eqref{1.1}$\thicksim$\eqref{1.3} gives
\beqq
\begin{aligned}
r^n\bigl(|a_n|+|b_n|\bigr)
&\leq
2
\left[
\frac{
\Gamma\left(\frac{2p-1}{2(p-1)}\right)
}{
\sqrt{\pi}\,
\Gamma\left(\frac{3p-2}{2(p-1)}\right)
}
\right]^{\frac{p-1}{p}}
M_p(r,f)\\
&\leq
2
\left[
\frac{
\Gamma\left(\frac{2p-1}{2(p-1)}\right)
}{
\sqrt{\pi}\,
\Gamma\left(\frac{3p-2}{2(p-1)}\right)
}
\right]^{\frac{p-1}{p}}
\|f\|_p.
\end{aligned}
\eeqq
Finally, letting $r\to1^-$ yields
\beqq
|a_n|+|b_n|
\leq
2
\left[
\frac{
\Gamma\left(\frac{2p-1}{2(p-1)}\right)
}{
\sqrt{\pi}\,
\Gamma\left(\frac{3p-2}{2(p-1)}\right)
}
\right]^{\frac{p-1}{p}}
\|f\|_p.
\eeqq

It remains to prove the sharpness part. We first consider
the case $p=1$. Fix $\xi\in\mathbb{T}$ and consider the Poisson kernel
\beqq
f_\xi(z)
=
\frac{1-|z|^2}{|\xi-z|^2}.
\eeqq
Since $f_\xi$ is positive and
$
M_1(r,f_\xi)=1
$ for all $r\in[0,1)$,
we have
\beqq
\|f_\xi\|_1=1.
\eeqq
Moreover,
\beqq
f_\xi(z)
=
1+
\sum_{j=1}^{\infty}\overline{\xi}^{\,j}z^j
+
\overline{
\sum_{j=1}^{\infty}\overline{\xi}^{\,j}z^j
}.
\eeqq
This gives
\beqq
a_j=b_j=\overline{\xi}^{\,j},
\qquad j\geq1,
\eeqq
and consequently
\beqq
|a_n|+|b_n|
=
2
=
2\|f_\xi\|_1.
\eeqq
Therefore, the constant $2$ in the case $p=1$ is sharp.

Finally, suppose that $1<p<\infty$. Fix $n\geq1$ and define
\beqq
\psi_{n,p}(t)
=
\left[
\frac{
\Gamma\left(\frac{2p-1}{2(p-1)}\right)
}{
\sqrt{\pi}\,
\Gamma\left(\frac{3p-2}{2(p-1)}\right)
}
\right]^{-1/p}
\operatorname{sgn}(\cos(nt))
|\cos(nt)|^{\frac{1}{p-1}},
\qquad 0\leq t\leq2\pi.
\eeqq
By \eqref{1.3}, we have
\beqq
\frac1{2\pi}
\int_0^{2\pi}
|\psi_{n,p}(t)|^p\,dt
=
\left[
\frac{
\Gamma\left(\frac{2p-1}{2(p-1)}\right)
}{
\sqrt{\pi}\,
\Gamma\left(\frac{3p-2}{2(p-1)}\right)
}
\right]^{-1}
\frac1{2\pi}
\int_0^{2\pi}
|\cos(nt)|^{\frac{p}{p-1}}\,dt=1.
\eeqq
Hence,
\beq\label{1.4}
\|\psi_{n,p}\|_{L^p(\mathbb{T})}=1.
\eeq

Let $f_{n,p}$ be the Poisson integral of $\psi_{n,p}$, namely,
\beqq
f_{n,p}(re^{i\theta})
=
\frac1{2\pi}
\int_0^{2\pi}
\frac{1-r^2}
{1-2r\cos(\theta-t)+r^2}
\psi_{n,p}(t)\,dt.
\eeqq
It follows from Theorem \ref{ThmB} and \eqref{1.4} that
\beqq
\|f_{n,p}\|_p=1.
\eeqq

Using the Fourier expansion of the Poisson kernel,
\beqq
%P_r(\theta-t) 
f_{e^{it}}(re^{i\theta})
=
1+
\sum_{k=1}^{\infty}
r^ke^{ik\theta}e^{-ikt}
+
\sum_{k=1}^{\infty}
r^ke^{-ik\theta}e^{ikt},
\eeqq
we obtain
\beqq
\begin{aligned}
f_{n,p}(re^{i\theta})
&=
\frac1{2\pi}
\int_0^{2\pi}\psi_{n,p}(t)\,dt+
\sum_{k=1}^{\infty}
\left(
\frac1{2\pi}
\int_0^{2\pi}
\psi_{n,p}(t)e^{-ikt}\,dt
\right)r^ke^{ik\theta}\\
&\quad+
\sum_{k=1}^{\infty}
\left(
\frac1{2\pi}
\int_0^{2\pi}
\psi_{n,p}(t)e^{ikt}\,dt
\right)r^ke^{-ik\theta}.
\end{aligned}
\eeqq
Comparing coefficients with the canonical expansion of $f_{n,p}$, we have
\beqq
a_k
=
\frac1{2\pi}
\int_0^{2\pi}
\psi_{n,p}(t)e^{-ikt}\,dt
\eeqq
and
\beqq
\overline{b_k}
=
\frac1{2\pi}
\int_0^{2\pi}
\psi_{n,p}(t)e^{ikt}\,dt.
\eeqq
In particular,
\begin{align}\label{1.5}
a_n
&=
\frac1{2\pi}
\int_0^{2\pi}
\psi_{n,p}(t)e^{-int}\,dt\\ \nonumber
&=
\left[
\frac{
\Gamma\left(\frac{2p-1}{2(p-1)}\right)
}{
\sqrt{\pi}\,
\Gamma\left(\frac{3p-2}{2(p-1)}\right)
}
\right]^{-1/p}
\frac1{2\pi}
\int_0^{2\pi}
\operatorname{sgn}(\cos(nt))
|\cos(nt)|^{\frac1{p-1}}
\bigl(\cos(nt)-i\sin(nt)\bigr)\,dt.
\end{align}
Since
\beqq
\frac{d}{dt}
|\cos(nt)|^{\frac{p}{p-1}}
=
-\frac{np}{p-1}
\operatorname{sgn}(\cos(nt))
|\cos(nt)|^{\frac1{p-1}}
\sin(nt),
\eeqq
we have
\begin{align}\label{1.6}
-i\int_0^{2\pi}
\operatorname{sgn}(\cos(nt))
|\cos(nt)|^{\frac1{p-1}}
\sin(nt)\,dt
&=
i\frac{(p-1)}{np}
\int_0^{2\pi}
\frac{d}{dt}
|\cos(nt)|^{\frac{p}{p-1}}\,dt\\
\nonumber
&=
i\frac{(p-1)}{np}
\left[
|\cos(nt)|^{\frac{p}{p-1}}
\right]\Big|_{0}^{2\pi}\\ \nonumber
&=0.
\end{align}
Consequently, \eqref{1.5} and \eqref{1.6} yield
\beqq
a_n
&=&
\left[
\frac{
\Gamma\left(\frac{2p-1}{2(p-1)}\right)
}{
\sqrt{\pi}\,
\Gamma\left(\frac{3p-2}{2(p-1)}\right)
}
\right]^{-1/p}
\frac1{2\pi}
\int_0^{2\pi}
|\cos(nt)|^{\frac{p}{p-1}}\,dt\\
&=&
\left[
\frac{
\Gamma\left(\frac{2p-1}{2(p-1)}\right)
}{
\sqrt{\pi}\,
\Gamma\left(\frac{3p-2}{2(p-1)}\right)
}
\right]^{\frac{p-1}{p}}.
\eeqq
A similar computation gives
\beqq
\overline{b_n}
=
\left[
\frac{
\Gamma\left(\frac{2p-1}{2(p-1)}\right)
}{
\sqrt{\pi}\,
\Gamma\left(\frac{3p-2}{2(p-1)}\right)
}
\right]^{\frac{p-1}{p}}.
\eeqq
Therefore,
\beqq
|a_n|+|b_n|
=
2
\left[
\frac{
\Gamma\left(\frac{2p-1}{2(p-1)}\right)
}{
\sqrt{\pi}\,
\Gamma\left(\frac{3p-2}{2(p-1)}\right)
}
\right]^{\frac{p-1}{p}}.
\eeqq
Since $\|f_{n,p}\|_p = 1$, equality is actually attained; therefore, the constant is sharp for each fixed $n \geq 1$.
The proof of this theorem is complete.
\qed

\section*{}

\bigskip

{\bf Data Availability} Our manuscript has no associated data.

{\bf Conflict of interest} The authors declare that they have no conflict of interest.

\bigskip

\section*{Acknowledgments}	The first author was partly supported by the
National Science Foundation of China (grant no. 12571080)
and Guangxi Natural Science Foundation $\#$2026GXNSFFA00640002.
The second author was partially supported by the National Natural Science Foundation of China (Grant No. 12371071) and the Key Project of the NSF of Hunan Province (Grant No. 2026JJ30002). The third author was partially supported by the Scientific Research Fund of the Hunan Provincial Education Department (Grant No. 25A0086). The fourth author was  also 
partially supported by the National Natural Science Foundation of China (Grant No. 12371071) and the Key Project of the NSF of Hunan Province (Grant No. 2026JJ30002).

\end{document}